# HISTORICAL ORIGINS OF THE NINE-POINT CONIC
# THE CONTRIBUTION OF EUGENIO BELTRAMI

**Maria Alessandra Vaccaro**

Università degli Studi di Palermo

**Abstract**

In this paper, we examine the evolution of a specific mathematical problem, i.e. the nine-point conic, a generalisation of the nine-point circle due to Steiner. We will follow this evolution from Steiner to the Neapolitan school (Trudi and Battaglini) and finally to the contribution of Beltrami that closed this journey, at least from a mathematical point of view (scholars of elementary geometry, in fact, will continue to resume the problem from the second half of the $19^{th}$ to the beginning of the $20^{th}$ century). We believe that such evolution may indicate the steady development of the mathematical methods from Euclidean metric to projective, and finally, with Beltrami, with the use of quadratic transformations. In this sense, the work of Beltrami appears similar to the recent (after the anticipations of Magnus and Steiner) results of Schiaparelli and Cremona. Moreover, Beltrami's methods are closely related to the study of birational transformations, which in the same period were becoming one of the main topics of algebraic geometry. Finally, our work emphasises the role played by the nine-point conic problem in the studies of young Beltrami who, under Cremona's guidance, was then developing his mathematical skills. To this end, we make considerable use of the unedited correspondence Beltrami – Cremona, preserved in the Istituto Mazziniano, Genoa.

**Introduction**

As is well known, Eugenio Beltrami was appointed in October 1862 to the chair in Complementary Algebra and Analytical Geometry of the University of Bologna even before obtaining his bachelor's degree. This decision was agreed upon by Francesco Brioschi and Luigi Cremona at a time when Beltrami had already shown his brilliant talents as a mathematician, but had produced very little.[1] As shown by the correspondence with Cremona,[2] kept in the archive of the Istituto Mazziniano in Genoa, Beltrami's mathematical apprenticeship did not end with the "conquest" of the chair, but continued for some years. Among the several topics with which he dealt in that period, we found it interesting to discuss a problem already studied by Jakob Steiner and later by Nicola Trudi: the nine-point conic and its connection with the geometry of the triangle and quadratic transformations.

---

[1] At that point of time Beltrami was the author of just three papers: (Beltrami 1861a, 102-108), (Beltrami 1861b, 257-283), (Beltrami 1861c, 283-284), the last one written as a letter to the editor of the journal, Prof. Barnaba Tortolini.

[2] About Luigi Cremona and his correspondence, see also http://www.luigi-cremona.it/.

Today this subject is widely studied from an educational point of view, also for its applications in dynamic geometry,[3] but its historical origins are often misrepresented. Indeed, as it often happens for elementary problems, the results of this kind of research have been rediscovered over time. For instance, Maxime Bôcher is often credited with the discovery,[4] in 1892, of the nine-point conic and its properties, thus postdating the history of this interesting configuration by about fifty years. Indeed, the nine-point conic was introduced in 1844 by Steiner in an article that appeared in Italian in a minor journal. This paper, which expounds properties stated but not proved by the Swiss mathematician, received some attention only in the Italian community. In 1856 and later in 1863, Trudi was the first to provide proofs for some of Steiner's theorems, followed by Giuseppe Battaglini and Giusto Bellavitis. The latter, a teacher in Padua, in his periodical Riviste di giornali scientifici (Reviews of scientific journals) 1862-1863, analysed Battaglini's work, also giving his own version of the problem of the nine-point conic.[5] In 1870, Pietro Cassani, teacher of Giuseppe Veronese at the technical institute of Venice, added to the subject, but it was especially Beltrami who between 1862 and 1874 not only fully proved the results by Steiner, but also generalised them and obtained complete results that were much deeper than those known up to then. So far we have traced in some detail the history of the nine-point conic to highlight the fact that, from the 1860s, the topic was well known to those who worked in elementary geometry in Italy, from Naples to Pavia and Padua. We claim that at the end of the 1850s and the beginning of the 1860s the study of problems in elementary geometry was very important for Italian mathematicians (mainly Cremona and Beltrami on the one side and Trudi and Battaglini on the other) to get fully acquainted with the new methods of projective geometry. The deep link between these kind of studies and the more advanced ones in algebraic geometry (Cremonian transformation, for example) have been studied in (Vaccaro 2016, 9-44).

Interest in the subject was resumed in the 1890s in the United States when, in a completely independent manner, Bôcher redefined the nine-point conic.[6] Several contributions followed[7] in which, again independently, some of the properties of the conic were rediscovered, culminating in a master's thesis[8] at Berkeley in 1912. In the same period, another series of papers on this subject appeared in Great Britain.[9] The first of these was due to Robert Edgar Allardice[10] who, in his article published in 1900, credited Beltrami for having first generalised the nine-point conic and stated that in subsequent works on this theme, written ignoring the contribution by the Italian mathematician, the same results were obtained again. A few years later Peter Pinkerton,[11] who seemed to ignore even the work of the British school, dealt with the same subject. From that moment on the topic was studied very frequently in the field of elementary geometry, but without reaching a deep vision such as Beltrami's.

---

3 See for instance (de Villiers, 2005) and (de Villiers, 2006).
4 See for instance (Pierce 2016, 27-78).
5 Besides Battaglini's paper, Bellavitis examined Beltrami's 1862 paper too. The two reviews appeared in Atti dell'Imperiale Regio Istituto Veneto di Scienze, Lettere ed Arti, (3), 8, 205-207 and 1281-1283.
6 Bôcher opened his paper as follows: "It does not seem to have been noticed that a few well-known facts, when properly stated, yield the following direct generalization of the famous nine-point circle theorem".
7 See for instance (Holgate 1893, 73-76) and (Gates 1894, 185-188).
8 The Master of Science thesis, The Nine-point Conic, was authored by Maud Minthorn under the supervision of Prof. M. W. Haskell and can be found at https://babel.hathitrust.org/cgi/pt?id=uc1.b3808276;view=1up;seq=11.
9 In fact, as early as 1858, Thomas Turner Wilkinson, a member of the Royal Astronomical Society, had extended the nine-point circle in an article published in The Lady's and Gentleman's Diary, in a similar way to what Steiner did. For more details see (Wilkinson 1858, 86-87).
10 (Allardice 1900, 23-32).
11 (Pinkerton 1905, 24-31).

The purpose of this paper is twofold: first of all, we want to highlight the influence of this topic on the Neapolitan geometry school and the role that this subject played in the early scientific training of Beltrami and, secondly, we will try to re-establish the real historical development of research on this topic. While Beltrami's later works on non-Euclidean geometry, differential geometry, and mathematical physics have been extensively studied by historians of mathematics, his studies and his youthful works are little known. Moreover, his correspondence with Jules Hoüel in (Boi, Giacardi and Tazzioli 1998) reveals his prolonged interest in this subject.[12] Through this work, therefore, we want to contribute to the historical studies on the mathematical approach of young Beltrami, also quoting extensively from his unpublished correspondence with Cremona. In conclusion, Beltrami's early work, now almost completely forgotten, is in our opinion rich in results that are much deeper than most of the subsequent rediscoveries on this subject, especially about the close connection between the nine-point conic and quadratic transformations.

**The nine-point circle**

The geometry of the triangle was born, as an independent discipline, in the second half of the nineteenth century and deals with problems concerning aligned or concyclic points with respect to a triangle. In this context, most of the main properties had already been highlighted by scholars of the emerging projective geometry, in particular by Steiner. One of the most interesting results in this direction concerns the so-called nine-point circle.

The nine-point circle, also known as "Euler circle" or "Feuerbach circle", after Karl Wilhelm Feuerbach, is constructed from an arbitrary triangle *ABC* and passes through nine special points: the midpoints of the sides, *D*, *E*, *F*, the feet of the altitudes, *H*, *I*, *J*, and the midpoints *K*, *L*, *M* of the segments that connect the vertices of the triangle to the orthocentre *O* (Fig. 1.).

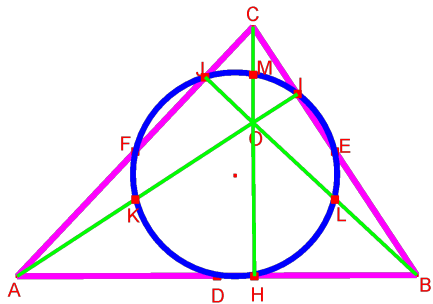

Fig. 1.

The earliest author to whom the discovery of the nine-point circle has been attributed is Leonhard Euler, but no one has been able to provide evidence or reference to any passage in Euler's writings in this regard. Based on historical research by John Sturgeon Mackay published in (Mackay 1892, 19-57), there have been several independent discoverers of the nine point circle. In (Feuerbach 1822) the theorem "the circle which goes through the feet of the perpendiculars of a triangle meets the sides at their mid points" is stated, but nothing is said of the other three points. In addition, Feuerbach proves that this circle touches the incircle and the excircles of the triangle (Fig. 2.). This last property is also stated in (Steiner

---
12 See also (Beltrami 1874, 543-566), (Beltrami 1876, 241-262), (Beltrami 1879, 233-312).

1828-29, 43). At the end of a long note in (Steiner 1833), Steiner claims that when he announced this theorem, he was not aware that it had been previously proved by Feuerbach.

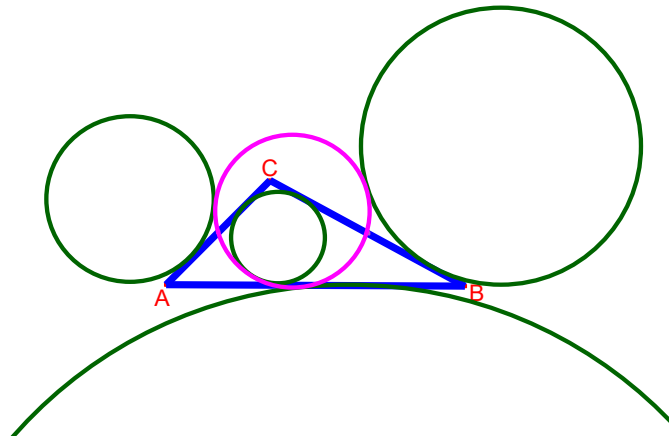

Fig. 2.

In his article (Steiner 1828-29, 37-64) Steiner shows, among other things, that the nine-point circle property is only a particular case of a more general theorem. The designation "nine-point circle" was given to it by Olry Terquem in (Terquem 1842, 198). From quotes in the relevant articles of the 1860s, it is clear that in 1862 a large set of results on the nine-point circle was well known to the Italian mathematicians, above all in the form used by Steiner and resumed in (Trudi 1863, 29-32), relating to the so-called complete quadrangle.

A complete quadrangle *ABCD* consists of the six lines that join two at a time the four vertices *A*, *B*, *C* and *D*. The intersection points $A_1$, $B_1$ and $C_1$ of the three pairs of opposite sides *AD* and *BC*, *AB* and *CD*, *AC* and *BD*, are the vertices of a triangle, called the fundamental triangle (Fig. 3.). In the case in which the pairs of opposite sides are orthogonal, i.e. if one of the vertices coincides with the orthocentre of the triangle formed by the remaining three, the complete quadrangle is said to be orthogonal (Fig. 4.).

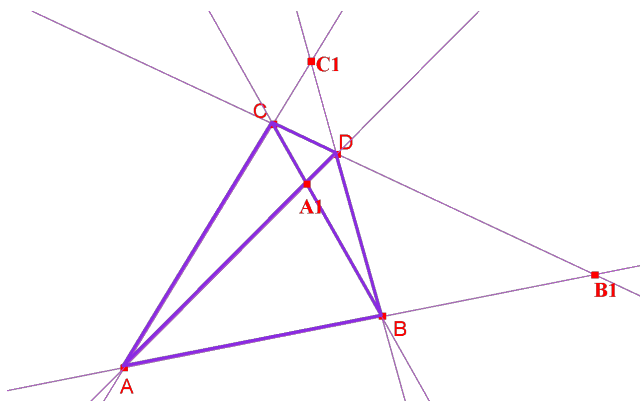 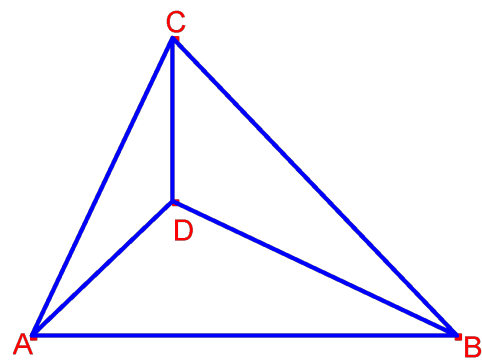

Fig. 3                              Fig. 4.

The midpoints of the diagonals *AD* and *BC* of the quadrangle and the midpoint of the segment $B_1C_1$, which connects the intersection points of the pairs of opposite sides, lie on the line *g*, nowadays called Gauss line of *ABCD* (Fig. 5.).

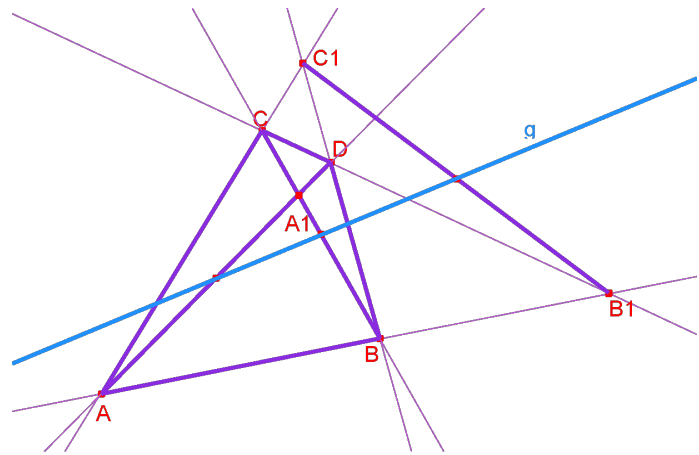

Fig. 5.

**Steiner's note**

Between September 1843 and the first half of 1844, Steiner, with Carl Borchardt, Peter Gustav Lejeune Dirichlet, Carl Jacobi, and Ludwig Schläfli, made a long journey (over nine months) in Italy, through Milan, Genoa, Pisa, Florence, Rome, Naples and other places. This journey was of great importance to Italian mathematicians, allowing the re-establishment of significant international contacts. The German and Swiss mathematicians came in contact, among others, in Pisa with Ottavio Mossotti and the physicist Carlo Matteucci (later both teachers of Enrico Betti), in Rome with Barnaba Tortolini (first editor of the Annals of Mathematics) and Domenico Chelini (later closely linked to Beltrami and Cremona in Bologna), in Naples with Vincenzo Flauti (teacher of Trudi), in Messina with Placido Tardy (one of the protagonists of Italian mathematics in the first period after Unification).[13]

In Rome, where they arrived after the Conference of Italian scientists in Lucca, both Jacobi and Steiner published works. In his article (Steiner 1844, 147-161), closely related to (Steiner 1828-29, 37-64) and published in Italian in a small journal,[14] Steiner dealt with the nine-point conic[15] studying two classical problems, dealing with the inscription and the circumscription about a given quadrilateral of a conic of maximum or minimum area. In particular, with regard to the question of the maximum-area conics inscribed in a quadrilateral, a problem in which Euler and Carl Friedrich Gauss had also been interested, Steiner stated that the centres of such conics belong to the so-called Gauss line. In addition, there are exactly two inscribed maximum-area conics, an ellipse and a hyperbola,[16] and the midpoint of their centres coincides with the centroid of the quadrilateral.

As regards the second problem, about minimum-area conics circumscribed about a quadrangle[17] *ABCD*, Steiner claims first that the locus of the centres of all the conics of the bundle through *A*, *B*, *C* and *D* is a conic $\Gamma$ passing through the midpoints of the six sides of the quadrangle and the points of intersection $A_1$, $B_1$ and $C_1$ of the three pairs of opposite sides (Fig. 6.).

---

[13] More information on this trip can be found in the biography of Jacobi in (Koenigsberger 1904).
[14] In a note in (Beltrami 1863a, 109), Beltrami claims that the memoir had been translated into Italian by Jacobi, who was helped by Prof. Domenico Chelini. This paper was also published in Italian, two years later, in the Journal of Crelle (Steiner 1846, 97-106).
[15] Steiner does not use the name "nine-point conic" that would be given, as we will see later, by Trudi and Battaglini only in 1863.
[16] In the case of the hyperbola, Steiner means the area of the corresponding ellipse with the same principal axes as the hyperbola.
[17] This term, due to Steiner, corresponds to German "vollständiges Viereck". The word Viereck, in English quadrangle, was translated into Italian in 1844 as "quadrigono".

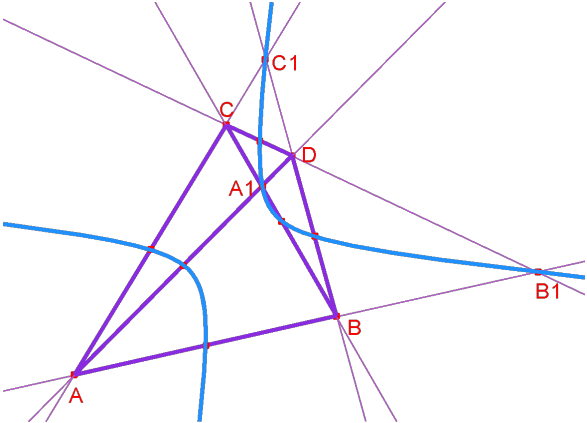

Fig. 6.

Next Steiner argues that finding the minimum-area conic among all the conics circumscribed about a quadrangle *ABCD* is equivalent to determining the minimum-area conic between all those circumscribed about the triangle *ABC* and having their centres on the conic $\Gamma$, a problem that he had already solved earlier in the same paper. Then the Swiss mathematician describes some properties of the conic $\Gamma$, including the fact that its centre coincides exactly with the centroid of the quadrangle.

Steiner's generalisation of the nine-point circle can be understood by looking again at the configuration of the nine-point circle in a slightly different way: rather than a triangle *ABC* and its orthocentre *D*, let us consider the complete orthogonal quadrangle *ABCD* (Fig. 7.). The definition of nine-point circle can be modified as follows: it is the locus of the centres of the conics circumscribed about an orthogonal quadrangle (all equilateral hyperbolas), which is the circle that passes through the midpoints of the six sides (*E*, *F*, *G*, *K*, *L*, *M*) and the three intersection points of opposite sides (*H*, *I*, *J*). Steiner would not note it explicitly, but, as we will see later, Trudi did.

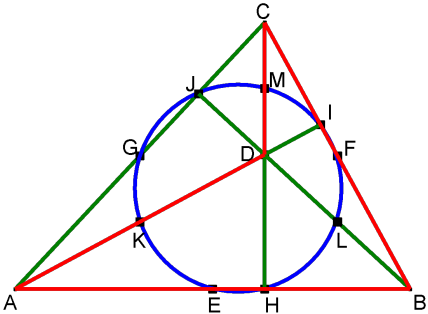

Fig. 7.

By removing the restriction for the quadrangle to be orthogonal, the theorem is still valid by replacing the word "circle" with "conic", that is, the nine-point conic, the locus of the centres of the bundle of conics through *A*, *B*, *C* and *D*. Steiner further notes that any conic $\Omega$ circumscribed about the fundamental triangle (the triangle $A_1B_1C_1$ in Fig. 8.) intersects each side of the complete quadrangle *ABCD* in a different point.

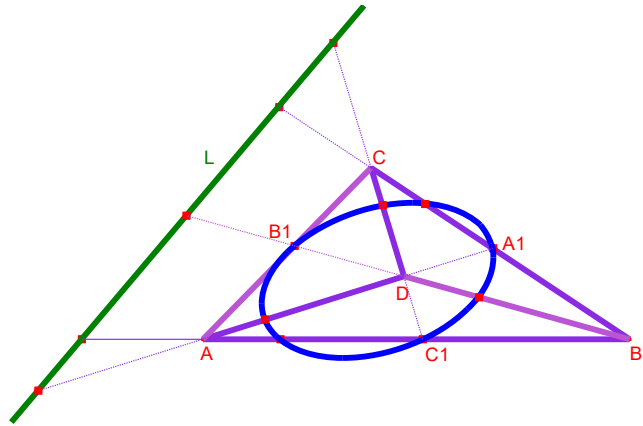

Fig. 8.

Furthermore, he considers, for each side, the harmonic conjugate point of the endpoints of the side and the intersection point. Then all six harmonic conjugates lie on the same line *L*: the poles of this line with respect to the conics of the bundle through *A*, *B*, *C* and *D* belong to the conic *Ω*. It is pointed out that when *Ω* is the nine-point conic, then the line *L* becomes the line at infinity. As we shall see later, this result will be a source of inspiration for Beltrami.

Steiner obtains other propositions from the previous ones by passing from the metric methods already used to projective processes.

Steiner generalises to the nine-point conic *Γ* the interesting property of the nine-point circle: given a triangle *ABC*, having called *D* its orthocentre, the nine-point circle common to the four triangles determined by points *A*, *B*, *C* and *D*, taken three at a time, is tangent to the sixteen circles inscribed in and escribed to the four triangles. Therefore, in each of the four triangles obtained by taking the vertices of the complete quadrangle three at a time, four conics similar to *Γ* may be inscribed and escribed and each of these sixteen conics is tangent to *Γ*.

An analytical proof of the fact that the locus of the centres of the conics circumscribed about a quadrangle is a conic is found in (Trudi 1856, 239-284). This proof became commonly known when Salmon inserted it as an exercise in (Salmon 1879, 153), the sixth edition of his treatise on conic sections. However, there is no reference to Steiner's statement, but, as already mentioned, between 1844 and 1892 the subject of the nine-point conic was widely studied.

**Trudi's contribution**

The first mathematician who showed an interest in Steiner's propositions and provided their proofs was Nicola Trudi,[18] who had met the Swiss mathematician, already in April 1844, on the occasion of his trip to Naples. Trudi claims, in the introduction of (Trudi 1856, 239), to have proved Steiner's statements quite early:

> Mr. Steiner made his singular solution known in a pamphlet published in Rome in 1844, which he presented to our friend Flauti, in the same year when he came with illustrious Iacobi [sic] to reside among us, having as its title: Teoremi su le coniche iscritte e circoscritte; but there, following his ordinary custom, he just announced the mere results; and, keeping silent about any analysis or proof, naturally impels the readers to find them out. Thus, driven to this work, of

---

18 About Nicola Trudi see (Ferraro, 2013).

> seventeen notable theorems stated therein, I could very soon assert the first sixteen of them.

His proofs were published only in 1856, even though they had already been presented to the Accademia in 1854. His first work on this subject is actually made up of two memoirs, the first devoted to the properties of the conics circumscribed about a quadrangle, the second to the question of the maximum or minimum area of the conics passing through four points. Trudi, after giving the definition of a complete quadrangle, distinguishes it from a complete quadrilateral.[19] He then defines a first-species quadrangle as one in which each vertex lies outside the triangle identified by the other three (Fig. 9a) and a second-species one[20] when only one vertex lies within the triangle determined by the remaining ones (Fig. 9b). He also denotes by three-intersection triangle, the triangle that has as its vertices the intersection points of the three pairs of opposite sides, that is, the fundamental triangle.

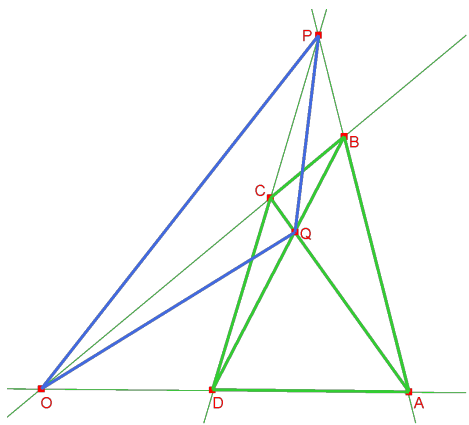 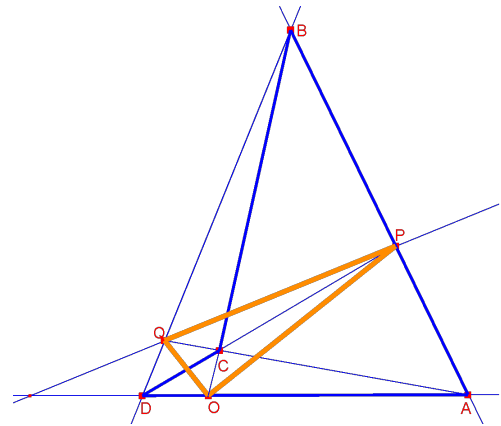

Fig. 9a.      Fig. 9b.

After recalling the property that in a complete quadrangle the midpoints of the three diagonals belong to the same straight line, the so-called Gauss line, Trudi derives that the midpoints of two opposite sides are aligned with the midpoint of the line segment joining the intersection points of the other two pairs of opposite sides (Figs. 10a and 10b). Using similar arguments, he deduces the following theorem: "In the complete quadrangle the three lines joining the midpoints of the three pairs of opposite sides meet in the same point, which divides each by half, and pass through the midpoints of the sides of the three-intersection triangle" (Trudi 1856, 244). This intersection point ($M$ in Figs. 10a and 10b) is the centroid of the four vertices of the quadrangle.

---

[19] Trudi denotes by "complete quadrilateral" the figure obtained by considering four straight lines, however lying in a plane, together with their intersection points.
[20] Complete orthogonal quadrangles belong to this second category.

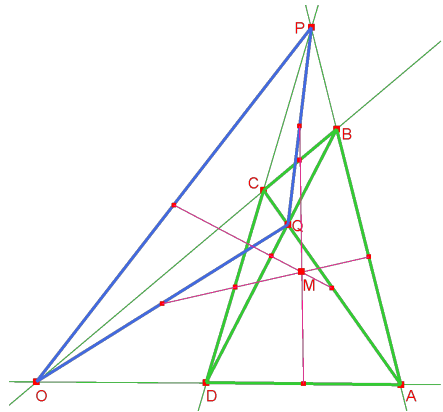 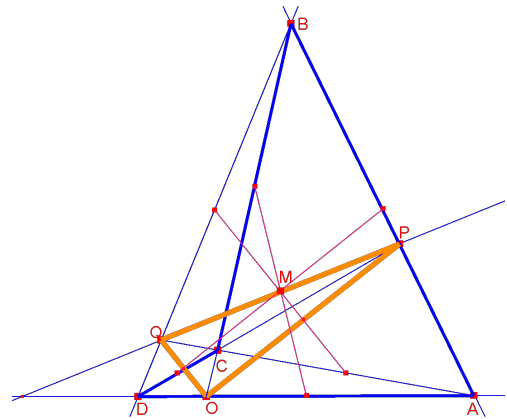

Fig. 10a.                    Fig. 10b.

Passing then to examine the type of conics that can be circumscribed about a quadrangle, he proves, using coordinate equations of the curves, that exactly two parabolas and infinitely many ellipses and hyperbolas pass through a first-species quadrangle, while infinitely many hyperbolas can be circumscribed about a second-species one. These hyperbolas are equilateral in the case of orthogonal quadrangles; more precisely, Trudi proves that the equilateral hyperbolas passing through the vertices of a triangle all meet in its orthocentre. Through analytical arguments, Trudi obtains the already known property: "in every quadrangle, the three-intersection triangle is such that each of its vertices is the pole of the corresponding opposite side with respect to any conic circumscribed about the quadrangle" (Trudi 1856, 253).

By studying the centre of a conic as the intersection of two of its diameters, he derives the equation of the locus of the centres of the conics passing through the vertices of a quadrangle and proves that "the locus of the centres of the innumerable conics circumscribable about a quadrangle is itself a conic section" (Trudi 1856, 253). In Trudi's opinion this locus deserves to be studied and discussed not just for its possible important applications, but also for its main properties, Trudi devotes the rest of this first memoir to examining the results that follow from them. First of all, he shows that this conic (locale), which he, like Steiner, does not name after the nine points, is circumscribed about the fundamental triangle and also passes through the midpoints of the sides of the quadrangle. Then he deduces that this locus is a hyperbola if the quadrangle is of the first species, an ellipse if the quadrangle is of second species: in particular, a circle if the quadrangle is orthogonal (Fig. 11b) and an equilateral hyperbola if the quadrangle is inscribed in a circumference (Fig. 11a.). In any case, it is a conic whose centre of symmetry coincides with the centroid of the quadrangle.

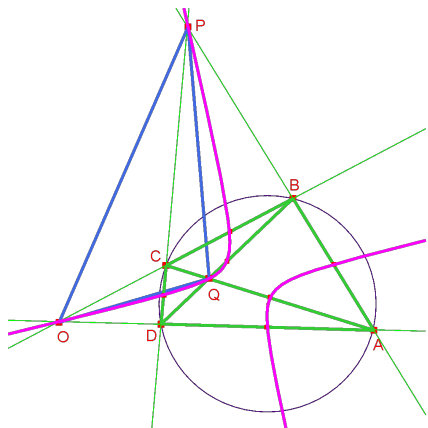 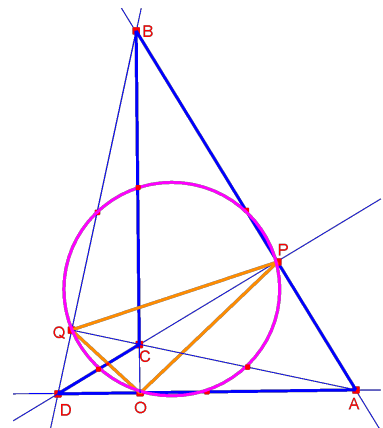

Fig. 11a.   Fig. 11b.

Trudi took up the subject again seven years later, in (Trudi 1863, 29-32), a paper devoted to the nine-point circle. In this article he provided a simple proof, based on trilinear coordinates, of the fact that the nine-point circle is tangent to each circle that is tangent to the three sides of the triangle. Thus, considering the figure consisting of the triangle and its altitudes as a complete orthogonal tetragon,[21] Trudi derives the fact that the nine-point circle relative to any of the four triangles into which the tetragon is decomposed is always the same; hence he states the following proposition:

> In the complete orthogonal tetragon the circle determined by the three intersections of the three pairs of opposite sides passes through the midpoints of all six sides. This circle, moreover, is the locus of the centres of all the conics (equilateral hyperbolas) circumscribed about the tetragon; and it is tangent to the sixteen inscribable circles, four at a time, in each of the four triangles determined by the vertices of the tetragon taken three at a time (Trudi 1863, 32).

Subsequently, by projectively transforming the configuration relative to the circle, he derives the theorem partly due to Jean-Victor Poncelet and partly to Steiner himself:

> The locus of the centres of the conics circumscribable about an arbitrary tetragon is another conic, which passes through the midpoints of all six sides and through the intersections of the three pairs of opposite sides. And this conic is tangent to the sixteen similar conics, placed in a similar way to it, which can be inscribed four at a time in each of the four triangles determined by the four vertices of the tetragon, taken three at a time (Trudi 1863, 32).

Trudi concludes by writing that "thus, the locus of the centres of the conics circumscribable about a tetragon can be called the nine-point conic, just like we say the nine-point circle" (Trudi 1863, 32).

In a slightly different context, in September 1862, Battaglini analysed the issue using almost the same words and proving that: "the locus of the centres of the conics circumscribed about a quadrangle is a conic that passes through the midpoints of its six sides, and through the intersection points of the three pairs of opposite sides. This conic will therefore be called the nine-point conic" (Battaglini 1862, 172). In his work Battaglini progressed quite far in the study of this conic, also proving the analogue of Steiner's theorem concerning the sixteen conics tangent to the nine-point conic. It is therefore to the Neapolitan community, Trudi and Battaglini, that we owe, besides the proof of the main properties, the name of this conic.

---

21 In footnote 1 to (Trudi 1863, 31-32), Trudi, after recalling the definitions of complete quadrilateral and quadrangle (quadrigono), admitting that both he and Battaglini had long used the latter name in several papers, claims that "recently, distinguished and learned Prof. Bellavitis has remarked that the Italian language has the word tetragono, more befitting than quadrigono to express the same idea; so he proposes the name of tetragono completo [= complete tetragon]".

**Beltrami's method**
In the introduction to (Beltrami 1862, 362), Beltrami made an explicit reference to Trudi's second paper, which "drew the scholars' attention to some very elegant theorems about what is called the nine-point circle". The theorems mentioned are those stated by Steiner in (Steiner 1828-29, 37-64) and developed in (Steiner 1844, 147-161). After also mentioning Battaglini among those who had dealt with this topic, Beltrami focussed on the fact that the proof provided by Trudi only applies to a particular case of the theorems stated by Steiner thus:

> In this short paper, I propose to deal again with this topic, considering it in its generality and taking the opportunity to examine very closely the nature and properties of a geometric transformation that is suggested as a spontaneous consequence of the previously proved theorems. The discussion of a special case of this transformation will show its intimate connection with other transformations that are well known and frequently used in geometry (Beltrami 1862, 362).

Thus Beltrami explicitly states the close connection existing between the nine-point conic and a particular quadratic transformation, which was called in (Hudson 1927) the standard quadratic transformation $T_2$.

After some preliminaries, he introduces the complete quadrangle *ABCD* and a line *r*. Let *E*, *F*, *G* be the intersection points of the opposite sides, forming the so-called fundamental triangle, $P_1$, $P_2$, $P_3$, $P_4$, $P_5$, $P_6$, the points where *r* intersects, respectively, the sides *AC*, *AB*, *BC*, *AD*, *CD*, *BD*, and $P'_1$, $P'_2$, $P'_3$, $P'_4$, $P'_5$, $P'_6$ the corresponding harmonic conjugates with respect to the previously mentioned sides (Fig. 12.). Having chosen as a projective reference system the points *E*, *F*, *G*, he determines the coordinates of the points and the equations of the geometrical entities of interest and proves that the points $P'_1$, $P'_2$, $P'_3$, $P'_4$, $P'_5$, $P'_6$, belong to the same conic: the nine-point conic corresponding to the straight line r.

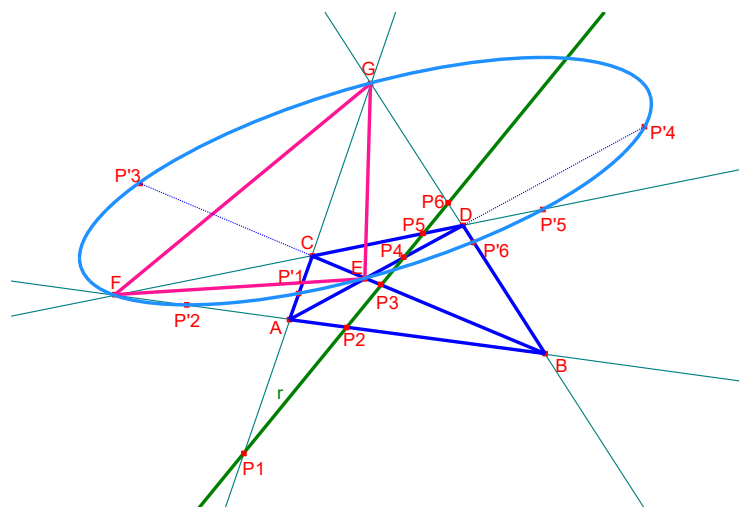

Fig. 12.

Thus, the theorem can be stated as follows:

> If in the plane of a complete quadrangle a transversal is drawn, and on each of its six sides the harmonic conjugate point of the point in which the side meets the transversal is determined, then the six points so determined lie on a conic, which also passes through the three intersection points of the opposite sides of the complete quadrangle (Beltrami 1862, 368).

Beltrami shows, again analytically, that the nine-point conic coincides with the locus of the poles of the line $r$ with respect to all the conics of the bundle determined by the vertices of a complete quadrangle. In this way he generalises the nine-point conic defined by Steiner, that is the locus of the centres of all the conics circumscribed about a complete quadrangle. Then in (Beltrami 1862, 361-395) he deduces Steiner's nine-point conic through a projective geometry argument: indeed, the centres of the conics of the bundle are just the poles of the line at infinity. From this theorem it can be deduced that the conics in the bundle that are tangent to $r$, which are known to be two, have their poles on the line $r$ itself, but also on the nine-point conic. Thus, the line $r$ intersects the nine-point conic precisely in the points $X$ and $X'$ where it is tangent to the two conics of the bundle that are tangent to $r$ (Fig. 13).

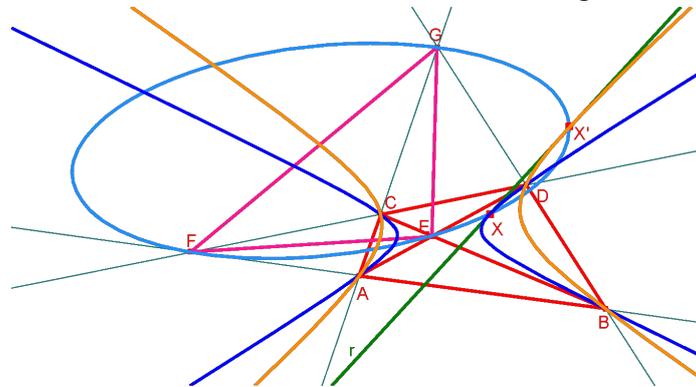

Fig. 13.

Moreover, since the pairs of points at which any line intersects the conics circumscribed about a quadrangle form a quadratic involution, Beltrami states the following theorem: "the double points of the involution that the conics circumscribed about a quadrangle determine on an arbitrary transversal are the points at which the latter meets the nine-point conic corresponding to it" (Beltrami 1862, 369).

Following this procedure, we can see how the situation turns out to be radically different from the one studied by Steiner: it is no longer a static configuration, but a dynamic one, in which, having fixed a quadrangle, it is possible to define a transformation which associates, to each point, a point and, to each line, a conic circumscribed about a fixed triangle. It is the classical context of quadratic transformations,[22] which precisely in that period were being subject to growing attention, and Beltrami would devote his attention exactly to making one of these transformations explicit. Thus, he changed the point of view by stating that: "it is important to remark that the previously proved properties could be simply deduced from the theorem that the polars of a point with respect to the infinitely many conics circumscribed about a quadrangle all pass through a same point" (Beltrami 1862, 370). Indeed, from this theorem, he deduced that, given the quadrangle $ABCD$, a line $r$ and a point $H$ on $r$, it is sufficient to determine the polars with respect to two conics of the bundle through $A$, $B$, $C$ and $D$, to obtain

---

[22] With regard to the historical aspects of the quadratic transformations, see (Hudson 1927) and (Snyder 1970).

*H′*, the corresponding of *H*. He chooses two degenerate conics, one given by the lines *AB* and *CD* and the other consisting of the lines *AC* and *BD* (Fig. 13). Beltrami shows that the locus of the points *H′*, when *H* moves on *r*, is precisely the nine-point conic relative to the line and the quadrangle. Indeed, the polar of an arbitrary point of *r* with respect to the degenerate conic (*AB*, *CD*) is a line passing through the intersection point *F* of the lines *AB* and *CD*, and this polarity is a projectivity between the line *r* and the sheaf through *F*. Similarly, by repeating the same argument with respect to the other degenerate conic (*AC*, *BD*), there is a projectivity between the line *r* and the sheaf through *G*, so the two sheaves of lines are in projective correspondence with each other. Thus, by using Newton's organic generation, the locus of the intersection points of lines that correspond to each other through a projectivity is a conic that we can call *Γ*. This conic, besides containing the points *F* and *G*, also passes through *E*, since we may repeat the same procedure also considering the degenerate conic consisting of the lines *AD* and *BC*, which intersect at *E*.

Arguing this way, having called *f* the correspondence that to each point *H* associates *H′* = *f*(*H*), that is, the point common to the polars of *H* with respect to all the conics of the bundle, we have a new way of generating the nine-point conic: the image of a line *r* under the transformation *f* defined above. Thus, Beltrami states:

> Given a complete quadrangle in a plane, each line of the same plane not in the quadrangle gives rise to a corresponding conic, circumscribed about the triangle formed by the intersection points of the three pairs of opposite sides of the quadrangle; and, reciprocally, each conic circumscribed about this triangle can be considered as corresponding to a unique and identified line of the plane (Beltrami 1862, 377-378).

Beltrami's point of view is totally new and much more modern than Steiner's. The transformation *f*, defined on the whole plane, which maps points to points and lines to conics circumscribed about a fixed triangle, is a quadratic transformation, and the nine-point conic fits into the framework of this theory.

Let us recall that in the early 1860s, on the momentum of the work by Giovanni Virginio Schiaparelli, the German mathematician Ludwig Immanuel Magnus (1790–1861) and Steiner himself, Luigi Cremona devoted himself to a profound revision of the concept of a quadratic transformation[23] that led him shortly after to the definition of birational transformations. It is therefore with satisfaction that Beltrami emphasised:

> Thus, we have here a correlation of points which proceeds according to this law, that to each point of the plane corresponds another unique and identified point of the same plane, and to each line corresponds a unique and identified conic circumscribed about a triangle that is invariable in its form and position, and conversely. This correlation is part of the more general one that was already discussed by several geometers, in particular by Steiner, Magnus and more recently the illustrious Mr. Prof. Schiaparelli (Beltrami 1862, 378-379).

In order to put Beltrami's work in its context, it is useful to recall that this March 1863 note was his fifth publication (the first one dating back to 1861), while he was already a professor at the University of Bologna. Moreover, Beltrami had already worked with Schiaparelli at the

---

[23] About this, see (Vaccaro 2016, 9-44).

Brera Observatory and was a close friend of Cremona, his colleague in Bologna, since his last years in Pavia as a student. Cremona published his famous memoir on birational transformations (Cremona 1863, 621-631) a few months after Beltrami, and in the same journal, and began by mentioning Magnus and Schiaparelli, just like Beltrami. Taking also into account the correspondence between the two great mathematicians mentioned in the next section, it can be assumed that it was precisely in the context of the conversations between them that ideas about this kind of transformation had matured.

Beltrami's paper continues by determining many interesting properties of the nine-point conic, among others, the sixteen conics tangent to that conic, which generalises the sixteen circles tangent to the nine-point circle:

> The sixteen conics that pass through the points common to an arbitrary line and the nine-point conic corresponding to it and that are inscribed in the four triangles formed by the six sides of the complete quadrangle are all touched by the nine-point conic (Beltrami 1862, 372).

Beltrami adds in a footnote: "This is the theorem that, for the special case of the orthogonal quadrangle and the transversal at an infinite distance, had been given by Feuerbach before Steiner" (Beltrami 1862, 372).

He also accurately determined the special cases in which the line is at infinity: the nine-point circle if the quadrangle is orthogonal, the Steiner conic if the quadrangle is arbitrary. Beltrami specifies that the only points fixed by the transformation are the four vertices of the quadrangle and that in each line there is only one pair of corresponding points: the intersection points of the straight line and the corresponding nine-point conic. Thus, "so that a line of the plane is touched by the corresponding conic, it is necessary (and sufficient) that it passes through one of the four vertices of the quadrangle, in which case the contact takes place in this same point" (Beltrami 1862, 380).

Beltrami proved a general theorem, which, as a special case, provides the definition of the transformation given previously, similar to that used today in algebraic geometry:

> If from a point of the plane the lines are drawn to the three intersection points of opposite sides of the quadrangle (the vertices of the fundamental triangle), and the harmonic conjugates of these lines are determined with respect to the sides passing through the respective intersection points, the three new lines thus obtained pass through one and the same point, which is the one corresponding to the given point. If the given point moves in the plane describing a line, the point determined in the aforementioned way describes a conic circumscribed about the fundamental triangle, and this conic is nothing but the conic corresponding to the line. All the points of a line passing through the intersection point of a pair of opposite sides of the quadrangle have their corresponding points in another line, passing through the same point and harmonically conjugate with the given line with respect to the two sides of the quadrangle (Beltrami 1862, 383).

This theorem is of considerable importance because it shows how Beltrami, starting from an elementary problem, detaches himself from the previous treatises, probably having in mind, like Cremona, applications to more advanced geometry topics.

Thus, by virtue of this theorem, having set $P' = f(P)$, the images of the points of the line $EP$ belong to the line $EP'$, which is none other than the symmetric of $EP$ with respect to any of

*AB* or *CD*, that is, the harmonic conjugate of the lines *AB*, *CD* (orthogonal to each other) and *EP*, as shown in Fig. 14.

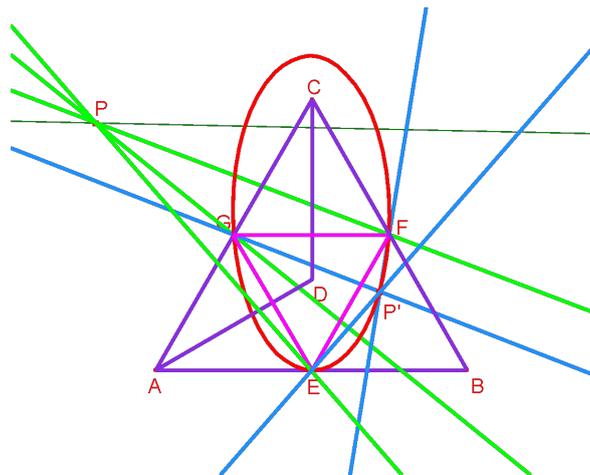

Fig. 14.

In the case when the quadrangle is orthogonal, the definition of this transformation coincides with that commonly used in 20th-century algebraic geometry texts,[24] that is, the image of a point *P* is determined as the intersection point of the symmetric lines of the lines that connect *P* with the vertices of the fundamental triangle with respect to the bisectors of the angles of the triangle itself. Moreover:

> If any one of the vertices of the fundamental triangle is joined by means of lines with any number of points on the plane, corresponding to one another two by two, a sheaf in involution is generated. The double lines of this involution are the two opposite sides of the quadrangle meeting at that point, and they are the lines that correspond to each other in this involution are the ones that pass through two corresponding points (Beltrami 1862, 384).

In other words, the lines fixed by *f* are those joining the vertices of the fundamental triangle with the four double points of the transformation itself, that is, the sides of the quadrangle. Beltrami explicitly states that the transformation *f* is biunivocal everywhere except in the points of the sides of the fundamental triangle; indeed, "to each point of one of the sides of the fundamental triangle the opposite vertex corresponds" and "to each vertex an arbitrary point of the opposite side corresponds" (Beltrami 1862, 384). Using modern terminology we could say that *f* blows up the vertices of the fundamental triangle and blows down the sides of the same triangle.

Beltrami points out that the image under *f* of a curve of order *n* has order $2n - (a + b + c)$, where *a*, *b* and *c* are the multiplicities with which the curve meets the vertices of the fundamental triangle. So a conic corresponds to a line or a conic or a cubic or a quartic, depending on whether the conic passes through three or two or a single vertex of the fundamental triangle or through none. Beltrami also remarks that the degree of the transformed curve is greater than or equal to $n/2$. Thus, to a cubic corresponds a conic if the curve passes through the three vertices and one of them is a double point, a cubic if the curve passes through the three vertices or through two vertices with a double point in one of them, a quartic if the curve passes through two vertices or if it has a double point in one of them, a quintic if the curve passes through a single vertex of the fundamental triangle, a sextic if it does not pass through any of the vertices.

---
24 See (Brieskorn and Knörrer 1986).

Beltrami's interest in this subject never ceased. Eleven years later, he returned to it. In the meantime Beltrami had turned his attention to differential geometry, producing important and well-known papers in non-Euclidean geometry (1868).

In the introduction to (Beltrami 1874, 543-544), he explains the reasons that led him to return several times to the same object of study, at different periods of his scientific activity. Beltrami had always tried to give a more "elegant" proof of the contact theorem, but the increasingly advanced developments in the methods of algebraic geometry had called on him to desist from this purpose. Let us remember that in this period Cremona, Alfred Clebsch and Max Noether gave a great boost to the study of birational transformations. However, as the subject continued to interest scholars, Beltrami decided to publish his research.

In this paper, the quadratic transformation we have mentioned becomes the starting point and, choosing as a projective reference system the one given by the orthogonal axes and the line to infinity, its analytical expression is:

$$f(x, y, z) = \left(\frac{1}{x}, \frac{1}{y}, \frac{1}{z}\right).$$

Hence, the fundamental triangle is $AOB$, where $A$ and $B$ are the points at infinity of the axes, while the quadrangle of the fixed points becomes $CDEF$, whose vertices are $(\pm 1, \pm 1, 1)$ (Fig. 15).

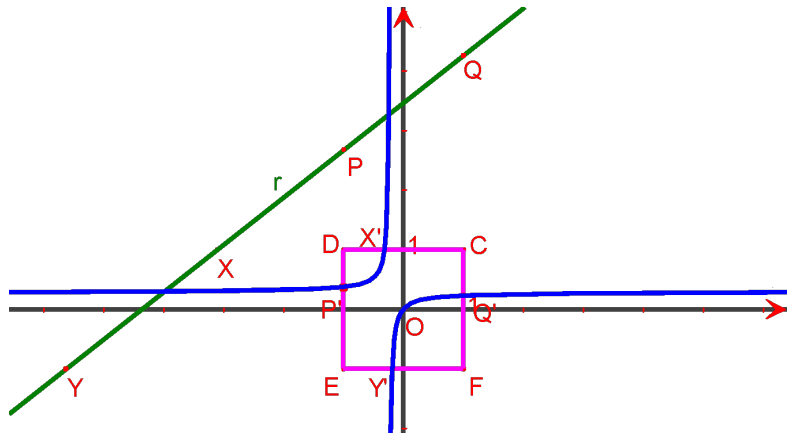

Fig. 15.

Denoting by $R \equiv ux + vy + wz = 0$ the equation of the line and with $\Gamma \equiv uyz + vxz + wxy = 0$ that of the corresponding nine-point conic, the bundle of conics determined by the quadrangle is represented by the equation: $px^2 + qy^2 + rz^2 = 0$, with $p + q + r = 0$. If the line $R$ passes through one of the vertices of the quadrangle, the conic passes through that point and is tangent to $R$ at it. If the line is the line at infinity, the nine-point conic is Steiner's and consists of just the axes. The image of the point $P$, in which the extension of the side $ED$ intersects $R$ (Fig. 15), under the transformation given by Beltrami is the point $P'$, which is the fourth point in harmonic ratio with $E$, $D$ and $P$.

Given the line $t$ that joins two corresponding points, $P$ and $P' = f(P)$, it is by construction the tangent to the conic in the bundle that passes through $P$. If $P$ is moved along a line, the lines $PP'$ have as their envelope a fourth-order, third-class curve $\Omega$, which is then projectively equivalent to the tricuspid hypocycloid, as stated by Beltrami: "it is but a projective transformation of the so-called tricuspid hypocycloid" (Beltrami 1874, 547) (Fig. 16).

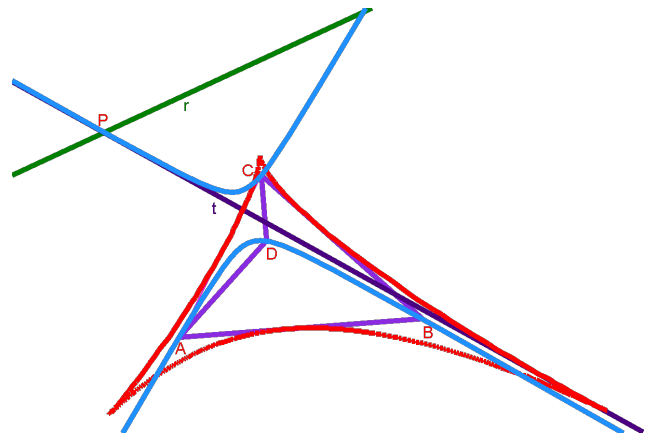

Fig. 16.

A correspondence is thus obtained between the projective plane and the dual plane, which associates with each line a quartic with three real or imaginary cusps. Beltrami studies in depth the quartic,[25] using precisely the transformation defined above. In particular, this curve admits a single bitangent (tangent in two distinct points), which is the line joining the two (real or imaginary) intersection points of the nine-point conic $\Gamma$ and the line corresponding to it under the quadratic transformation $f$.[26] By virtue of Plücker's formulas, Beltrami notes that, since the quartic has a bitangent, it has as its dual curve a rational cubic with a double point. He remarks that many of the properties of the quartic can be derived from this duality, as, for example, that the curve $\Omega$ cannot have any other singularity beyond the three cuspidal points. Beltrami also describes another construction of the quartic $\Omega$ by generating it pointwise rather than as an envelope: denoting by $P''$ the second intersection point of the line $PP'$ and the nine-point conic $\Gamma$ (the first being obviously $P'$), the locus of the points $X$, the harmonic conjugate points of $P''$ with respect to $P$ and $P'$, is a quartic equivalent to the tricuspid hypocycloid (Fig. 17).

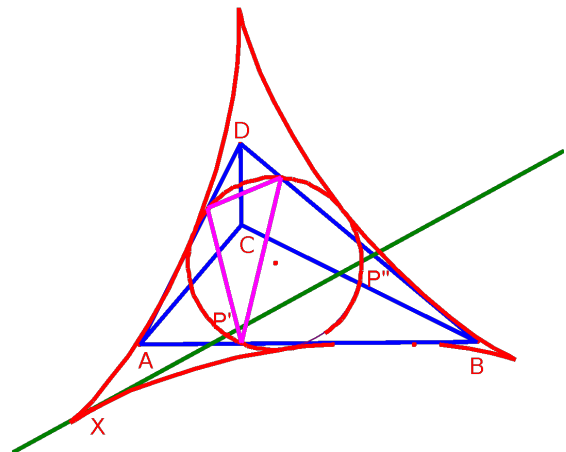

Fig. 17.

As can be seen, Beltrami's vision has expanded: the nine-point conic as the corresponding of a line through a quadratic transformation is no longer the subject of investigation, but the general properties of this quadratic transformation are analysed and, through them, higher-order curves are studied.

---

25 For this question, see (Palladino and Vaccaro 2018, 61-92).
26 In the case of the tricuspid hypocycloid the bitangent is imaginary, but in the more general case in which the line is proper, Beltrami's statement can be clearly visualised.

**Beltrami's ideas through his correspondence with Cremona**

As already mentioned in the introduction, in October 1862 Beltrami was appointed, before graduating, as a professor at the University of Bologna, and here he met again his friend and colleague Cremona. The correspondence[27] between the two mathematicians clearly shows not only the intimate friendship between them, but also Cremona's role as Beltrami's mentor in those years. The former is indeed generous with advice, encouragement and recommendations to his friend. As can be seen from his letter, dated Bologna 31 December 1860, Cremona, who used purely synthetic methods, did not want to induce his friend, who instead preferred analytical methods, to follow his approach (see below). At the same time, he discouraged Beltrami from writing a textbook that would not benefit his career:

> Since you left me complete freedom to act, I wrote Brioschi (on the 25th) so that he may find a way for you to graduate as soon as possible and with the least possible expense. But since he has come to Milan, my letter will not yet have reached him. Let us await his answer. [...] As for you, make sure of graduating soon. In fact, perhaps the appointment can be achieved even without a degree: but, in the uncertainty, it is good to obtain every useful qualification. I can only approve your study plan. I am not a man who discriminates against alternative ideas: as long as you show <u>firmness</u> and <u>resolve</u>, follow the inspiration of your intelligence: it is the best guide. Study the different branches of analysis, master the powerful tools created by modern science: then you will choose the field that best suits the nature of your mind. I am under the impression that you are intending to write an educational work about determinants – this I cannot approve. Mind this: study for yourself, do not worry about others now. The idea of making a good book, of simplifying and making a difficult-looking theory accessible to others, is seductive for a young man: I had that idea too, but fortunately Brioschi dissuaded me and I am grateful to him. Study for the sole purpose of educating <u>only yourself</u>: if you will have the opportunity to write something, let them be short papers and containing new ideas: papers that, without robbing you of much time, bring you much honour. When you have completely mastered the science, then you will occupy yourself as you see best.[28]

---

27 The unedited correspondence Beltrami – Cremona is preserved in the archive of the Istituto Mazziniano in Genoa (IMG). A special thank-you for the transcription of the correspondence to Dr. Paola Testi Saltini.

28 Avendomi tu lasciato piena libertà d'agire, ho scritto (il 25) a Brioschi perché trovi modo di laurearti il più presto e colla minor spesa possibile. Ma essendo egli venuto a Milano, la mia lettera non gli sarà ancor pervenuta. Attendiamo la risposta. […] Tu pensa a laurearti presto. Giacché, forse si potrà conseguire la nomina anche senza laurea: ma nel dubbio è bene premunirsi d'ogni titolo utile. Io non posso che approvare il tuo piano di studi. Io non sono uomo di idee esclusive: purché ci metta <u>fermezza</u> e <u>costanza</u>, segui pure l'inspirazione del tuo ingegno: essa è la miglior guida. Studia i diversi rami d'analisi, fa tuoi i poderosi strumenti creati dalla scienza moderna: poi sceglierai il campo che più si confarà all'indole della tua mente. Mi pare di aver capito che tu sia nell'intenzione di fare un lavoro didattico sui determinanti - questo non so approvare. Bada bene: studia per te, non curarti ora degli altri. L'idea di fare un buon libro, di semplificare e rendere accessibile altrui una teoria d'aspetto difficile è seducente per un giovane: io pure ebbi tale idea, ma fortunatamente Brioschi me ne distolse ed gliene son grato. Studia coll'unico scopo d'erudire <u>te solo</u>: se ti si offrirà occasione di fare alcun lavoro, siano lavori brevi e contenenti idee nuove: lavori che senza rubarti molto tempo ti procaccino molto

In the letter sent from Bologna on 11 November 1861, Cremona congratulates Beltrami on his first paper (Beltrami 1861a, 102-108), which would be published shortly:

> On the same day when I sent the aforementioned letter, Confalonieri arrived from Piacenza, and he gave me your envelope containing your paper. I have examined it and I have had it also examined by Chelini: we both congratulate you, since it is a work carried out with great skill. We believe you can send it to Tortolini to be included in the Annali. In view of this publication, we have taken the liberty of making some modifications that do not affect the substance: modifications that you will keep or not, as you prefer.[29]

On the other hand, Beltrami exchanged with Cremona ideas and points of view regarding his research, as can be seen from the letter excerpts, dated Milan 16 March 1862, that are given below. With the sentence "I have but little familiarity with those powerful methods of geometric procreation",[30] Beltrami refers to the methods of the emerging theory of quadratic transformations developed in those years by Cremona and Schiaparelli:

> After having written my previous letter, a generalisation of the theorem you know about presented itself to me too. And this generalisation, although it is less broad than what you are stating, has this peculiarity, that it does not fit into your general theorem. I obtained it with a very brief argument in analytical geometry, which I ask you to verify, lest I were in error. [...] It is therefore probable that this second theorem might give rise to a new generalisation that includes in itself that which you have given. I do not dare to do it, because I have but little familiarity with those powerful methods of geometric procreation (allow me this phrase), of which you have spoken to me so many times. Still, I hope I have not fallen into error in establishing the theorem. If it holds, however, it could be interesting in the study of a transformation method in which the point of coordinates a : b : c in one figure is made to correspond with that of coordinates 1/a : 1/b : 1/c in another, a transformation that would be a special case of that indicated by Magnus, in which to each point of one figure one point in the other corresponds, and consequently, to all lines of one figure, conics circumscribed about the same triangle in the other.[31]

---

onore. Quando sarai completamente padrone della scienza, allora ti occuperai come meglio ti piacerà (IGM, 051-11298 (7986)).

29 Nello stesso giorno in cui impostai la lettera menzionata, arrivò Confalonieri da Piacenza, e mi consegnò un tuo piego contenente il tuo lavoro. L'ho esaminato io e l'ho fatto esaminare anche dal Chelini: entrambi te ne facciamo le nostre congratulazioni, poiché è un lavoro condotto con molta abilità. Crediamo che tu possa mandarlo al Tortolini per essere inserito negli Annali. In vista di tale inserzione, ci siamo presi la libertà di fare alcune modificazioni che non intaccano però la sostanza della cosa: modificazioni però che tu rispetterai o no, come credi (IGM, 051-11723 (8411)).

30 The term procreation here is used as a synonym for generation.

31 Dopo averti scritto la precedente mia lettera, si presentò a me pure una generalizzazione del teorema che sai. E questa generalizzazione, benché sia meno ampia di quella che tu enunci, ha questo di particolare che non rientra punto nel tuo teorema generale. Io l'ho ottenuta con brevissime considerazioni di geometria analitica, che ti prego di verificare, per vedere se mai fossi in errore. [...] È probabile quindi che questo secondo teorema possa dar luogo ad una nuova generalizzazione che comprenda in sé quella da te indicata. Io non mi arrischio a farla, perché ho poca famigliarità con quei potenti metodi di procreazione geometrica (permettimi la frase) dei quali tu mi hai parlato tante volte. D'altronde non vorrei esser caduto in errore nello stabilire il teorema. Se questo poi sussistesse, esso potrebbe essere interessante nello studio di un metodo di trasformazione nel quale al punto di coordinate a : b : c in una figura, si facesse corrispondere quello di coordinate 1/a : 1/b : 1/c nell'altra,

Cremona's reply from Bologna is dated 3 April 1862:

> I sent your proof to Bellavitis. The beautiful theorem that you announced to me in your last letter is correct, and from it we can deduce the following [...]. You should devote yourself to this beautiful question, which must certainly lead to remarkable results. I do not have time to think about it.[32]

In the letter that follows, dated 22 February 1863, Beltrami, who was already a temporary professor in Bologna, recounted to his interlocutor:

> Today I am sending Battaglini two small papers of mine for publication in Naples Giornale. You already know one of them, about the visual angle of a second-order surface,[33] which I have not seen published in Terquem['s journal] and which I believe was lost. You partly know the other one too, and it is the proof of the main theorems about the nine-point conics, which I have completed by adding the part about the contacts of these conics with the 16 ones you know. What gave me an opportunity to return for a moment on this topic was a note published in this regard by Trudi,[34] in which he expresses a wish that those theorems be given a simple proof, and begins by giving one himself. This proof of his is simple enough indeed, but it seems to me that he took the thing backwards. He proves analytically the property for the circle and then extends it to conics by using purely geometrical considerations (parallel projections). This method would tend to imply that the analytical proof of this second part is more difficult than the first one's, while the opposite is true. In analytic geometry, the simplicity of calculations derives rather from a suitable choice of the axes and the coordinates than from special hypotheses made on the nature of the curves. Trudi's method seems to me to be more appropriate to research in pure geometry, where if an argument can be simplified it is by simplifying the objects of the argument itself, especially since pure geometry possesses its own, most fruitful methods of generalisation. But it seems to me that when we resort to analysis we must obtain all necessary generality with it alone, and geometry must intervene only to guide the operations and to set out the equations. This, it seems to me, is the true combination of analysis and geometry; the other one is nothing but a mixture, as the chemists say.
>
> Moreover, in going back to these topics, I have found that they can be extended to space, and that there, instead of a nine-point conic, we have a third-order surface on thirty-two points and nine lines. This will form the subject of another paper,[35] unless the publication of this first article is dragged out, in which case I will redo the whole thing together. But, as it almost always happens when moving from

---

trasformazione che sarebbe un caso particolare di quella indicata da Magnus, in cui a ciascun punto di una figura corrisponda un punto nell'altra, e per conseguenza, a tutte le rette dell'una figura, coniche circoscritte ad uno stesso triangolo nell'altra (IGM, 049-10232 (6929)).

32 Ho mandato la tua dimostrazione al Bellavitis. Il bel teorema che tu mi comunichi nell'ultima tua è esatto e da esso si può dedurre quanto segue [...]. Tu dovresti occuparti di questa bella quistione, la quale deve condurre certamente a rimarchevoli risultati. Io non ho tempo di pensarci (IGM, 049-10233 (6930)).

33 This paper was (Beltrami 1863c, 68-73).

34 Beltrami refers to (Trudi 1863, 29-32).

35 In May 1863, a few months after this letter, Beltrami published the paper (Beltrami 1863b, 208-217 and 354-360).

plane to space, the ways of generalising those theorems can be multiple. The one I have found so far, however, seems to be the most spontaneous, for the perfect analogy of the analytical results too, although some properties are lost and some theorems do not find corresponding ones. Be that as it may, if I get something together about it, I'll let you know it before publishing it in any way.[36]

It can be seen that three years after the first letter, Beltrami's reflection on methods has matured. He now expresses in a profound way his view on the relationship between analytical and synthetic methods in geometry, a discussion that we believe is still very topical.

**Conclusion**

As already mentioned in the introduction, Steiner's 1844 article and those by the Italian authors on the nine-point conic did not obtain a great response and were largely rediscovered later, especially in the British cultural sphere.[37] Beltrami's 1862 note aroused instead interest in the German community, mainly leading to the development of a few meaningful results for the geometry of the triangle, but collateral with respect to the work of the Italian geometer. According to Gino Loria's remarks published in (Loria 1901, 392-440), the results obtained by Beltrami that aroused most interest in Germany were basically two:

I) the centre of the circle circumscribed about a triangle is the centroid of the centres of the circles inscribed in [and escribed to] the same triangle;

II) if from the vertices of a triangle three parallel lines are drawn and for each of them the symmetrical line with respect to the corresponding bisector is taken, three lines are found that are concurrent in one point.

---

36 Oggi stesso spedisco a Battaglini due piccoli miei lavori perché sian pubblicati nel Giornale di Napoli. Uno è quello che già conosci, relativo all'angolo visuale di una superficie di 2° ordine[36], che non veggo comparire nel Terquem e che credo andato smarrito. L'altro ti è pure in parte noto, ed è la dimostrazione dei principali teoremi relativi alle coniche dei nove punti, che ho completato coll'aggiungere la parte che riguarda i contatti di queste coniche colle 16 che sai. Mi diede occasione a ritornare per un momento sopra questo argomento una nota pubblicata in proposito dal Trudi[36], nella quale esterna il desiderio che di quei teoremi si dia una dimostrazione semplice, ed incomincia col darne una lui medesimo. Questa sua dimostrazione è abbastanza semplice per verità, ma mi pare che egli abbia preso la cosa a rovescio. Egli dimostra analiticamente la proprietà pel circolo e poi la estende alle coniche valendosi di considerazioni puramente geometriche (projezioni parallele). Questo metodo tenderebbe a far credere che la dimostrazione analitica di questa seconda parte sia più difficile di quella della prima, mentre è piuttosto vero il contrario. Nella geometria analitica la semplicità dei calcoli deriva piuttosto dalla scelta opportuna degli assi e delle coordinate che da ipotesi speciali fatte sulla natura delle curve. Il metodo del Trudi mi parrebbe opportuno per le ricerche di pura geometria, in cui se può semplificarsi il ragionamento è col semplificare gli oggetti del ragionamento stesso, tanto più che la geometria pura possiede poi i metodi di generalizzazione suoi propri e fecondissimi. Ma mi pare che quando si ricorre all'analisi si debba ottenere con essa sola tutte le necessarie generalità, e la geometria non debba venire che a guidare le operazioni e ad intavolare le equazioni. Questa, sembra a me, è la vera combinazione dell'analisi e della geometria; l'altra non ne è che una mescolanza, come dicono i chimici.

Nel ripigliare poi queste considerazioni ho trovato ch'esse possono estendersi allo spazio, e che ivi, invece di una conica di nove punti, abbiamo una superficie del terzo ordine di trentadue punti e di nove rette. Ciò formerà argomento di un altro scritto[36], a meno che la pubblicazione di questo primo articolo non vada per le lunghe, nel qual caso rifonderò insieme il tutto. Ma come avvien quasi sempre quando si passa dal piano allo spazio, le maniere di generalizzare quei teoremi pajono essere molteplici. Quella che ho trovato finora sembra però la più spontanea, anche per la perfetta analogia dei risultati analitici, benché alcune proprietà si perdano ed alcuni teoremi non trovino i loro riscontri. Comunque sia, se arriverò a metter insieme qualche cosa in proposito, la farò conoscere a te prima di pubblicarla in un modo qualsiasi (IGM, 049-09782 (6478)).

37 Apparently, the English-language authors, whom we have mentioned at the beginning, independently coined the "nine-point conic" designation.

The first theorem was proved analytically in 1864 by Johann August Grunert;[38] the proof was followed by a series of papers all published in the same year 1865 and in the same journal (by Grunert): Archiv der Mathematik und Physik.[39] For the second theorem an analytical proof was also provided by Grunert,[40] to whom Beltrami, who had obtained the result using methods in projective geometry, replied with a letter to the editor[41] proving its correctness through simple arguments in Euclidean geometry. We find further developments related to the geometry of the triangle, but without explicit reference to Beltrami's work, in the French school since the second half of the 1860s.[42]

On the other hand, with regard to the quadratic transformation, which plays a fundamental role in Beltrami's approach, there was no follow-up except in Cremona, who names it Beltrami's transformation in (Cremona 1865, 88-91).

In 1863, Beltrami devoted himself to extending this topic to three-dimensional space (Beltrami 1863b, 208-217 and 354-360). The extension to space was resumed independently in 1863 by Eugène Prouhet, who merely extended Feuerbach's theorem by determining the twelve-point sphere, and later in 1872 by Friedrich Emil Eckhardt.[43] The subject of the extension to space would occupy Beltrami's interest for a long time, as evidenced by the correspondence with Cremona.

Although the articles we have analysed show that Beltrami had made most of Cremona's methods his own, in fact his interest in this problem kept alive only at the level of specific "elementary" problems. On the other hand, Beltrami's modest approach in his 1874 paper, which was subtitled analytical exercise for good reasons, might have contributed to being underestimated by the main European mathematicians. Perhaps one could hypothesize that these papers were too general to interest scholars of triangle geometry and too elementary for mathematicians now interested in the classification of curves, surfaces and hypersurfaces in hyperspaces.

**References**

Allardice, R.E., 1900. On the nine-point conic. Proceedings of the Edinburgh Mathematical Society, 19, 23-32.

Battaglini, G., 1862. Nota sopra alcune questioni di geometria. Rendiconti della Reale Accademia delle Scienze Fisiche e Matematiche di Napoli, 1, 168-178.

Beltrami, E., 1861a. Intorno ad alcuni sistemi di curve piane. Annali di matematica pura e applicata, serie 1, tomo IV, 102-108.

Beltrami, E., 1861b. Sulla teoria delle sviluppoidi e delle sviluppanti. Annali di matematica pura e applicata, serie 1, tomo IV, 257-283.

Beltrami, E., 1861c. Di alcune formole relative alla curvatura della superficie. Annali di

---

38 (Grunert 1864, 354-356).
39 E.J. Noeggerath 89-91, R. Lobatto 234-235, C. Schmidt 238-241, C.G. Reuschle 364 and 483-487, vol. 43. C. Struve 119, C. Schmidt 120-124, W. Stammer 335-336, vol. 44.
40 (Grunert 1865, 102-108).
41 (Beltrami 1865, 481-483).
42 See (Romera-Lebret 2014, 253-302).
43 See (Prouhet 1863, 132-138) and (Eckhardt 1872, 30-49).

matematica pura e applicata, serie 1, tomo IV, 283-284.

Beltrami, E., 1862. Intorno alle coniche dei nove punti e ad alcune quistioni che ne dipendono. Memorie dell'Accademia delle Scienze dell'Istituto di Bologna, serie II, vol. II, 361-395.

Beltrami, E., 1863a. Sulle coniche di nove punti. Giornale di Matematiche, 1, 109-118.

Beltrami, E., 1863b. Estensione allo spazio a tre dimensioni dei teoremi relativi alle coniche dei nove punti. Giornale di Matematiche, 1, 208-217 and 354-360.

Beltrami, E., 1863c. Soluzione d'un problema relativo alle superficie di second'ordine. Giornale di Matematiche, 1, 68-73.

Beltrami, E., 1865. Auszug aus einem Briefe des Herrn Professor Eugenio Beltrami in Pisa an den Herausgeber, betreffend die im Archiv mitgetheilten Beltrami'schen Sätze. Archiv der Mathematik und Physik, 43, 481-483.

Beltrami, E., 1874. Intorno ad alcuni teoremi di Feuerbach e di Steiner. Esercitazione analitica. Memorie dell'Accademia delle Scienze dell'Istituto di Bologna, serie III, tomo V, 543-566.

Beltrami, E., 1876. Considerazioni analitiche sopra una proposizione di Steiner. Memorie dell'Accademia delle Scienze dell'Istituto di Bologna, serie III, tomo VII, 241-262.

Beltrami, E., 1879. Ricerche di Geometria analitica. Memorie dell'Accademia delle Scienze dell'Istituto di Bologna, serie III, tomo X, 233-312.

Bôcher, M., 1892. On a Nine Point Conic. Annals of mathematics, 6, 132.

Boi, L., Giacardi, L., Tazzioli, R., 1998. La découverte de la géométrie non euclidienne sur la pseudosphère, Les lettres d'Eugenio Beltrami à Jules Hoüel (1868-1881). Éditions Albert Blanchard, Paris.

Brieskorn, E., Knörrer, H., 1986. Plane algebraic curves. Birkhäuser Verlag, Basel.

Cremona, L., 1863. Sulle trasformazioni geometriche delle figure piane. Memorie dell'Accademia delle Scienze dell'Istituto di Bologna, 2, 621-631.

Cremona, L., 1865. On normals to conics, a new treatment of the subject. The Oxford, Cambridge, and Dublin Messenger of Mathematics, 3, n. X, 88-91.

de Villiers, M., 2005. A generalization of the nine-point circle and Euler line. Pythagoras, 62, 31-35. DOI: 10.4102/pythagoras.v0i62.112

de Villiers, M., 2006. The nine-point conic: a rediscovery and proof by computer. International Journal of Mathematical Education in Science and Technology, 37, n. 1, 7-14. DOI: 10.1080/00207390500138025

Eckhardt, F.E., 1872. Beiträge zur analytischen Geometrie des Raumes, Mathematische Annalen, 5, 30-49.

Ferraro, G., 2013. Excellens in arte non debet mori, Nicola Trudi da napolitano a italiano. https://halshs.archives-ouvertes.fr/halshs-00682088, 1-16.

Feuerbach, K.W., 1822. Eigenschaften einiger merkwürdigen Punkte des geradlinigen Dreiecks und mehrerer durch sie bestimmten Linien und Figuren. Eine analytisch-trigonometrische Abhandlung (Monograph ed.), Nürnberg.

Gates, F., 1894. Some considerations on the nine point conic and its reciprocal. Annals of Mathematics, 8, 185-188.

Grunert, J.A., 1864. Satz vom ebenen Dreieck. Archiv der Mathematik und Physik, 42, 354-356.

Grunert, J.A., 1865. Über einen geometrischen Satz vom Dreieck. Archiv der Mathematik und Physik, 43, 102-108.

Holgate, T., 1893. On the Cone of Second Order Which is Analogous to the Nine-Point Conic. Annals of Mathematics, 7, 73-76.

Hudson, H.P., 1927. Cremona transformations in plane and space. Cambridge University Press, Cambridge.

Koenigsberger, L., 1904. Carl Gustav Jacob Jacobi: Festschrift zur Feier der Hundertsten Wiederkehr seines Geburtstages. Druck und Verlag von B. G. Teubner, Leipzig.

Loria, G., 1901. Eugenio Beltrami e le sue opere matematiche. Bibliotheca Mathematica, 2, Issue 3, 392-440.

Mackay, J.S., 1892. History of the Nine-point Circle. Proceedings of the Edinburgh Mathematical Society, 11, 19-57.

Palladino, N., Vaccaro, M.A., 2018. L'ipocicloide tricuspide: il duplice approccio di Luigi Cremona ed Eugenio Beltrami. Bollettino di Storia delle Scienze Matematiche, 38, Issue 1, 61-92. DOI: 10.19272/201809201003

Pierce, D., 2016. Thales and the nine-point conic. The De Morgan Gazette, 8, n. 4, 27-78.

Pinkerton, P., 1905. On a nine point conic. Proceedings of the Edinburgh Mathematical Society, 24, 24-31.

Prouhet, E., 1863. Analogies du triangle et du tetraèdre. Cercle des neuf points, sphère des douze points. Nouvelles Annales de mathématiques, Série 2, Tome 2, 132-138.

Romera-Lebret, P., 2014. La nouvelle géométrie du triangle à la fin du XIX$^e$ siècle: des revues mathématiques intermédiaires aux ouvrages d'enseignement. Revue d'histoire des mathématiques, 20, 253-302. DOI: 10.24033/rhm.183

Salmon, G., 1879. A treatise on conic sections: containing an account of some of the most important modern algebraic and geometric methods. Longmans, Green and co., London.

Snyder, V., 1970. Selected topics in algebraic geometry. AMS Chelsea Publishing Company, New York.

Steiner, J., 1828-1829. Géométrie pure. Développement d'une série de théorèmes relatifs aux sections coniques. Annales de Mathématiques pures et appliquées, 19, 37-64.

Steiner, J., 1833. Die geometrischen Constructionen, ausgeführt mittelst der geraden Linie und eines festen Kreises. Ferdinand Dümmler, Berlin.

Steiner, J., 1844. Teoremi relativi alle coniche iscritte e circoscritte. Giornale Arcadico di Scienze Lettere ed Arti, 295, 147-161.

Steiner, J., 1846. Teoremi relativi alle coniche iscritte e circoscritte. Journal für die reine und angewandte Mathematik, 30, 97-106.

Terquem, O., 1842. Considérations sur le triangle rectiligne. Nouvelles annales de mathématiques, Série 1, Tome 1, 196-200.

Trudi, N., 1856. Memorie relative alle proprietà delle curve del 2° ordine circoscrittibili ad un quadrigono ed alla ricerca della minima tra esse in superficie. Memorie della Reale

Accademia delle Scienze, Napoli, 239-284.

Trudi, N., 1863. Nota intorno ad una proprietà del cerchio de' nove punti. Giornale di Matematiche, 1, 29-32.

Vaccaro, M.A., 2016. Dalle trasformazioni quadratiche alle trasformazioni birazionali. Un percorso attraverso la corrispondenza di Luigi Cremona. Bollettino di Storia delle Scienze Matematiche, 36, Issue 1, 9-44.

Wilkinson, T.T., 1858. Notae Geometricae. The Lady's and Gentleman's Diary, 86-87.